\documentclass[12pt,a4paper]{article}
\usepackage{amsmath,amsthm,amssymb,amscd}
\usepackage{graphicx}

 \DeclareMathOperator{\id}{id}
 
 \DeclareMathOperator{\sgn}{sgn}

\newcommand{\C}{{\mathbb{C}}}
\newcommand{\Z}{{\mathbb{Z}}}
\newcommand{\M}{{\mathcal{M}}}
\newcommand{\pair}[2]{\langle #1,#2\rangle}
\newcommand{\pairrho}[2]{\langle #1,#2\rangle_\rho}

\newtheorem{thm}{Theorem}[section]
\newtheorem{lem}[thm]{Lemma}

\theoremstyle{definition}
\newtheorem{defn}[thm]{Definition}
\theoremstyle{remark}
\newtheorem{rem}[thm]{Remark}

\begin{document}

\title{Faithfulness of the Lawrence representation of braid groups}
\author{Hao Zheng}
\date{}
\maketitle

\begin{abstract}
The Lawrence representation $L_{n,m}$ is a family of homological
representation of the braid group $B_n$, which specializes to the
reduced Burau and the Lawrence-Krammer representation when $m$ is
$1$ and $2$. In this article we show that the Lawrence
representation is faithful for $m \geq 2$.
\end{abstract}

\section{Introduction}

In \cite{Bigelow} and \cite{Krammer}, Bigelow and Krammer proved
via different approaches that the Lawrence-Krammer representation
of braid groups is faithful thus the braid groups are linear. In
fact, the Lawrence-Krammer representation is the only known
faithful representation of the braid group $B_n$ for $n\geq4$ till
now.

In this article, by making use of a reflexive representation
recently found by the author (ref. \cite{Zheng}), we generalize
the faithfulness of the Lawrence-Krammer representation to its
full family, the Lawrence representation (ref. \cite{Lawrence}).

\begin{thm}\label{thm:faith}
The Lawrence representation if faithful for $m \geq 2$.
\end{thm}

In the article, the Lawrence representation is defined
alternatively as follows. Let $B_n$ denote the Artin's $n$-strand
braid group (ref. \cite{Birman}), with standard generators $\{
\sigma_1,\dots,\sigma_{n-1} \}$, and set
$$B_{n,m} = \langle \sigma_1,\dots,\sigma_{n-1}, \sigma_n^2,
 \sigma_{n+1},\dots,\sigma_{n+m-1} \rangle \subset B_{n+m}.$$
They are the fundamental groups of
\begin{eqnarray*}
 && X_n = \{ (x_1,\dots,x_n) \mid x_i \in \C,
 x_i \neq x_j, \forall i \neq j \} / \Sigma_n, \\
 && X_{n,m} = \{ (x_1,\dots,x_{n+m}) \mid x_i \in \C, x_i \neq x_j,
 \forall i \neq j \} / \Sigma_n \times \Sigma_m
\end{eqnarray*}
respectively, where $\Sigma_n$ denotes the symmetric group of $n$
symbols.

Let $\xi_{n,m}$ be the reflexive representation over a free $\Z
B_{n,m}$-module $M_{n,m}$ defined in \cite{Zheng} (see Section
\ref{sec:def}). Let $q,t \in \C$ be two algebraically independent
numbers and let
$$\rho_{n,m} : \Z B_{n,m} \to \C$$
denote the ring homomorphism given by
$$\left\{ \begin{array} {lll}
\sigma_1,\dots,\sigma_{n-1} & \mapsto & 1, \\
\sigma_n^2 & \mapsto & q, \\
\sigma_{n+1},\dots,\sigma_{n+m-1} & \mapsto & t.
\end{array}\right.$$
The {\em Lawrence representation} is defined as the representation
$$L_{n,m} = \rho_{n,m} \circ \xi_{n,m}$$
over the $\C$-linear space
$$M^L_{n,m} = \C \otimes_{\rho_{n,m}} M_{n,m}.$$

\begin{rem}
It was shown in \cite{Zheng} that $L_{n,2}$ is precisely the
Lawrence-Krammer representation and it is easily derived from the
explicit matrix elements calculated in \cite{Zheng} that $L_{n,1}$
is precisely the reduced Burau representation (ref.
\cite{Birman}).
\end{rem}

\begin{rem}
It is known that the reduced Burau representation is faithful for
$n\leq3$ and not faithful for $n\geq5$ (ref. \cite{Bigelow1}), but
the case $n=4$ still remains open. Therefore, Theorem
\ref{thm:faith} shows that the faithfulness of Lawrence
representation is only unclear for $L_{4,1}$.
\end{rem}

Our proof essentially follows Bigelow's approach. In Section
\ref{sec:def}, we give a quick review of the reflexive
representation $\xi_{n,m}$. In Section \ref{sec:pairing}, we
define the pairing of noodles with multiforks and relate it to the
Lawrence representation via the notion of linear function. It is
the crucial part of the article. In Section \ref{sec:proof}, after
some preliminary lemmas prepared, the main theorem is established.

\section{A quick review of the representation $\xi_{n,m}$}\label{sec:def}

Let $D$ be a $2$-disk and $P=\{p_1,\dots,p_n\} \subset D \setminus
\partial D$ be a set of $n$ punctures. The space
$$Y_{n,m} = \{ (y_1,\dots,y_m) \mid y_i \in D \setminus P, \;
 y_i \neq y_j, \; \forall i \neq j \} / \Sigma_m$$
is homotopy equivalent to the fiber of the fiber bundle $X_{n,m}
\to X_n$, whose fundamental group is
$$\langle A_{1,n+1},\dots,A_{n,n+1}, \sigma_{n+1},\dots,\sigma_{n+m-1} \rangle
 \subset B_{n,m}$$
where $A_{i,j}$ is the standard pure braid defined by
$$A_{i,j} = \sigma_{j-1} \cdots \sigma_{i+1} \sigma_i^2
 \sigma_{i+1}^{-1} \cdots \sigma_{j-1}^{-1}.$$

Recall that an equivalent definition of $B_n$ is the mapping class
group $\M(D,P;\partial D)$, the group of all orientation
preserving homeomorphism $h : D \to D$ such that $h(P)=P$ and
$h|_{\partial D}=\id$, modulo isotopy relative to $P \cup \partial
D$. Regarding $\M(D,P;\partial D)$ and $\pi_1(Y_{n,m})$ as
subgroups of $B_{n,m}$ in the standard way, we have
$$\beta_*(\alpha) = \beta^{-1}\alpha\beta, \;\;\;
 \forall \beta \in \M(D,P;\partial D),\;\; \alpha \in \pi_1(Y_{n,m}).$$

\bigskip

\begin{figure}[h]
    \centering
    \includegraphics[scale=.6]{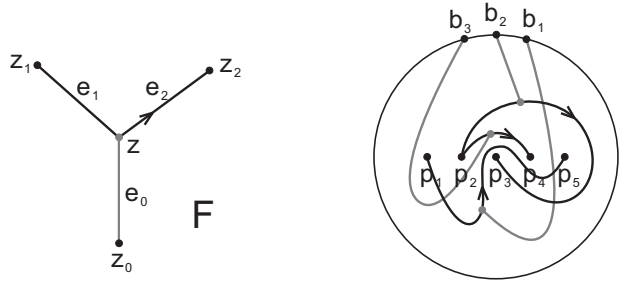}
    \caption{Complex $F$ and a multifork.}
    \label{fig:fig21}
\end{figure}

Let $F$ be the $1$-complex shown in Fig. \ref{fig:fig21}. It
consists of four $0$-cells $\{z,z_0,z_1,z_2\}$ and three $1$-cells
$\{e_0,e_1,e_2\}$. Let $e_t = e_1 \cup z \cup e_2$ and $e_h = z_0
\cup e_0 \cup z$ denote the {\em tine edge} and the {\em handle}
of $F$.

\begin{defn}
A {\em fork} is a map $\phi : F \to D$ such that $\phi|_{e_t}$ is
an embedding, $\phi(F) \cap \partial D = \phi(z_0)$ and $\phi(F)
\cap P = \{\phi(z_1),\phi(z_2)\}$. A {\em multifork} with $m$
components is an $m$-tuple of forks $\Phi=(\phi_1,\dots,\phi_m)$
such that both $\phi_1(e_t),\dots,\phi_m(e_t)$ are disjoint and
$\phi_1(e_h),\dots,\phi_m(e_h)$ are disjoint.
\end{defn}

\begin{defn}
Two forks $\phi$ and $\psi$ are called {\em homotopic}, denoted by
$\phi \simeq \psi$, if there is a homotopy $h_t : F \to D$ such
that $h_0=\phi$, $h_1=\psi$, $h_t(z_0)$ is independent of $t$, and
$h_t$ is a fork for all $0 \leq t \leq 1$. Two multiforks $\Phi =
(\phi_1,\dots,\phi_m)$ and $\Psi = (\psi_1,\dots,\psi_m)$ are
called {\em homotopic}, also denoted by $\Phi \simeq \Psi$, if
there are fork homotopies $h_{k,t} : \phi_k \simeq \psi_k$ such
that $(h_{1,t},\dots,h_{m,t})$ is a multifork for all $0 \leq t
\leq 1$.
\end{defn}

\medskip

Choose a base point $[b_1,\dots,b_m] \in Y_{n,m}$ where
$b_1,\dots,b_m \in \partial D$. Set
$$\Gamma_{n,m} = \{ (\phi_1,\dots,\phi_m) \mid \phi_i(z_0)=b_i, \;
 1 \leq i \leq m \}$$
and denote by $M^0_{n,m}$ the free $\Z B_{n,m}$-module generated
by $\Gamma_{n,m}$. Define four relations on $M^0_{n,m}$ as
follows.

\medskip

$R_H$: $\Phi_1 \sim \Phi_2$ if they are homotopic.

$R_R$: $(\phi_1,\dots,\phi_k,\dots,\phi_m) \sim
-(\phi_1,\dots,\phi_k r,\dots,\phi_m)$ where $r : F \to F$ denotes
the cell isomorphism that swaps $e_1$ and $e_2$.

\begin{figure}[h]
    \centering
    \includegraphics[scale=.6]{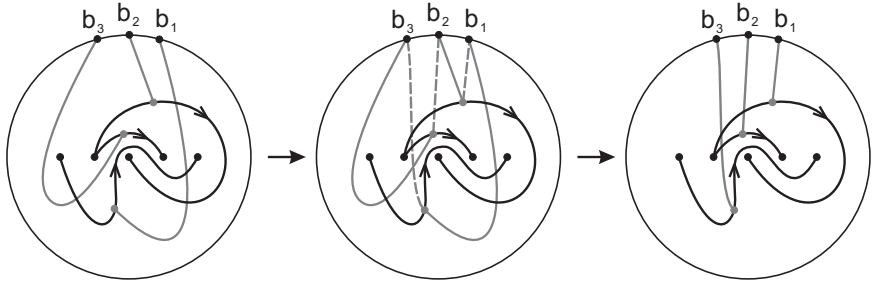}
    \caption{Relation $R_T$.}
    \label{fig:fig21}
\end{figure}

$R_T$: $(\phi_1,\dots,\phi_m) \sim \sgn\eta \cdot \alpha \cdot
(\varphi_1,\dots,\varphi_m)$ where $\eta \in \Sigma_m$ and $\alpha
\in \pi_1(Y_{n,m})$ if $\phi_k|_{e_t}=\varphi_{\eta(k)}|_{e_t}$
for all $1 \leq k \leq m$ and $\alpha$ is represented by the loop
that runs from $[b_1,\dots,b_m]$ to $[\phi_1(z),\dots,\phi_m(z)]$
along the curve $\{[\phi_1(t),\dots,\phi_m(t)] \mid t \in e_h\}$
and backs to $[b_1,\dots,b_m]$ along the curve
$\{[\varphi_1(t),\dots,\varphi_m(t)] \mid t \in e_h\}$.

\begin{figure}[h]
    \centering
    \includegraphics[scale=.6]{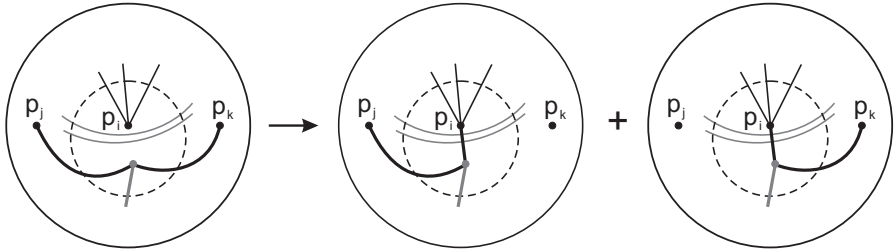}
    \caption{Relation $R_S$.}
    \label{fig:fig23}
\end{figure}

$R_S$: $\Phi \sim \Phi_1 + \Phi_2$ if $\Phi$ can be {\em split}
into $\Phi_1$ and $\Phi_2$ by doing a surgery on the tine edge of
a fork as shown in Fig. \ref{fig:fig23}.

\bigskip

Now set $M_{n,m} = M^0_{n,m} / (R_H,R_R,R_T,R_S)$. It turns out
that the action of $B_n$ on $M_{n,m}$ gives rise to a
representation over a finitely generated free $\Z B_{n,m}$-module.

\begin{thm}
$M_{n,m}$ is a finitely generated free $\Z B_{n,m}$-module.
Moreover, the action
$$\xi_{n,m}(\beta) : [\Phi] \mapsto [\beta \cdot \beta(\Phi)], \;\;\;
 \forall \Phi \in \Gamma_{n,m}, \;\; \beta \in B_n$$
gives rise to a representation of $B_n$ over $M_{n,m}$.
\end{thm}

\section{Pairing and linear function}\label{sec:pairing}

\begin{defn}
A {\em noodle} is an embedded oriented arc $N \subset D \setminus
P$ such that $\partial D \cap N = \partial N$ and all the points
$b_1,\dots,b_m$ lies to its left.
\end{defn}

\begin{figure}[h]
    \centering
    \includegraphics[scale=.6]{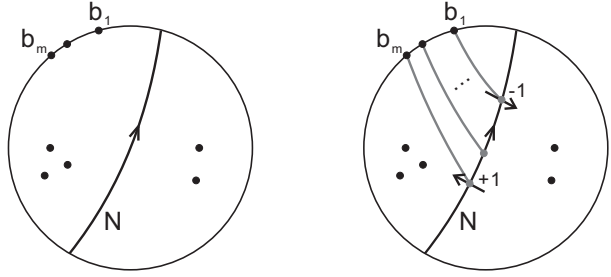}
    \caption{Noodle $N$, signs of intersections and a set of disjoint arcs.}
    \label{fig:faith21}
\end{figure}

\begin{defn}\label{defn:pairing}
Let $N$ be a noodle and $\Phi = (\phi_1,\dots,\phi_m)$ be a
multifork such that the tine edge of $\phi_i$ intersects $N$
transversely at $\{ x_{i,1},x_{i,2},\dots,x_{i,l_i} \}$. The {\em
pairing} of $N$ with $\Phi$ is defined as
$$\pair{N}{\Phi} =
 \sum_{i_1=1}^{l_1} \cdots \sum_{i_m=1}^{l_m}
 \epsilon_{1,i_1} \cdots \epsilon_{m,i_m}
 \alpha_{i_1,\dots,i_m}
 \in \Z \pi_1(Y_{n,m}),
$$
where $\epsilon_{j,i_j}$ is the sign of the intersection
$x_{j,i_j}$ of $N$ with $\phi_j(e_t)$, $\alpha_{i_1,\dots,i_m} \in
\pi_1(Y_{n,m})$ is represented by the loop that runs from
$[b_1,\dots,b_m]$ to $[\phi_1(z),\dots,\phi_m(z)]$ along the
handles of $\Phi$ (i.e. $\{[\phi_1(t),\dots,\phi_m(t)] \mid t \in
e_h\}$), then to $[x_{1,i_1},\dots,x_{m,i_m}]$ along the tine
edges of $\Phi$ (the subarcs of $\phi_j(e_t)$ from $\phi_j(z)$ to
$x_{j,i_j}$), and backs to $[b_1,\dots,b_m]$ along the disjoint
arcs shown in Figure \ref{fig:faith21}.
\end{defn}

It is straightforward to verify that, via the pairing, each noodle
$N$ gives rise to a $\Z B_{n,m}$-linear function
$$\pair{N}{\,\cdot\,} : M_{n,m} \to \Z B_{n,m},$$
and, further, a $\C$-linear function
$$\pairrho{N}{\,\cdot\,} : M^L_{n,m} \to \C.$$

Note that we have
$$\pairrho{N}{L_{n,m}(\beta) \cdot [\Phi]}
  = \pairrho{N}{[\beta(\Phi)]}, \;\;
  \forall \beta \in B_n, \;
  \Phi \in \Gamma_{n,m}.
$$
Especially, if $\beta$ is an element of the kernel of the Lawrence
representation $L_{n,m}$,
$$\pairrho{N}{[\Phi]}
  = \pairrho{N}{[\beta(\Phi)]}, \;\;
  \forall \Phi \in \Gamma_{n,m}.
$$

\medskip

\begin{rem}
For $m=2$, the last equation is precisely a generalization of
\cite[Basis Lemma]{Bigelow}. Here we obtain the equation via the
language of representation, which makes the topological meaning
much more accessible.
\end{rem}

\section{Proof of faithfulness}\label{sec:proof}

\begin{figure}[h]
    \centering
    \includegraphics[scale=.6]{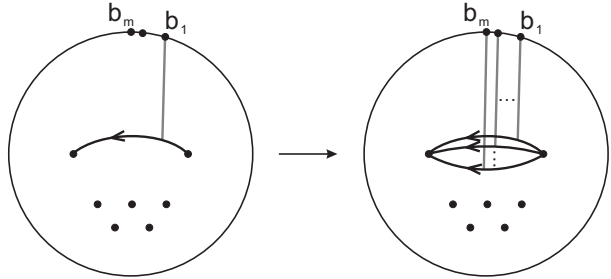}
    \caption{Fork to multifork.}
    \label{fig:faith31}
\end{figure}

In this section, let all forks $\phi$ satisfy $\phi(z_0) = b_1$
and denote by $\phi^{(m)}$ the multifork constructed from $m$
parallel copies of $\phi$ as shown in Figure \ref{fig:faith31}.

\begin{lem}\label{lem:faith:coef}
Let $N$ be a noodle and $\phi$ be a fork. Suppose the tine edge of
$\phi$ intersects $N$ transversely at $l$ distinct points and
\begin{eqnarray*}
 && \pair{N}{\phi^{(m)}} =
 \sum_{i_1,\dots,i_m=1}^{l}
 \epsilon_{i_1} \epsilon_{i_2} \cdots \epsilon_{i_m}
 \alpha_{i_1,\dots,i_m}, \\
 && \rho_{n,m}(\alpha_{i_1,\dots,i_m}) =
 q^{a_{i_1,\dots,i_m}}
 (-t)^{b_{i_1,\dots,i_m}},
\end{eqnarray*}
where $\epsilon_i$ and $\alpha_{i_1,\dots,i_m}$ are same as
Definition \ref{defn:pairing}. Then we have
\begin{eqnarray*}
 && \epsilon_i = (-1)^{b_{i,i}}, \\
 && a_{i_1,\dots,i_m} = \sum_{j=1}^m a_{i_j}, \\
 && b_{i_1,\dots,i_m} = \sum_{1 \leq j <k \leq m} b_{i_j,i_k}.
\end{eqnarray*}
\end{lem}
\begin{proof}
Note that for $j>k$, $b_{i_j,i_k}$ is the crossing number (define
the crossing number of the generator $\sigma_i^{\pm1}$ to be
$\pm1$) between the $(n+j)$-th and the $(n+k)$-th strand of the
braid $\alpha_{i_1,\dots,i_m}$, $a_{i_j}$ is the linking number
(half of the crossing number) of the $(n+j)$-th strand with the
former $n$ strands of the braid $\alpha_{i_1,\dots,i_m}$.

The identities follow from the facts that the crossing number
between the $(n+1)$-th and the $(n+2)$-th strand of $\alpha_{i,i}$
is even if and only if $\epsilon_i$ is positive,
$a_{i_1,\dots,i_m}$ is the linking number of the last $m$ strands
with the former $n$ strands of $\alpha_{i_1,\dots,i_m}$,
$b_{i_1,\dots,i_m}$ is the sum of the pairwise crossing numbers of
the last $m$ strands of $\alpha_{i_1,\dots,i_m}$, respectively.
\end{proof}

\begin{lem}\label{lem:faith:key}
Let $N$ be a noodle, $\phi$ be a fork and $m \geq 2$ be an
integer. If $\pairrho{N}{\phi^{(m)}} = 0$ then the tine edge of
$\phi$ is isotopic to relative to $\partial D \cup P$ to an arc
which is disjoint from $N$.
\end{lem}
\begin{proof}
Applying a preliminary isotopy, we may assume that the tine edge
of $\phi$ intersects $N$ transversely at $l$ distinct points where
$l$ is minimal in possible. Suppose $l>0$. In the notation of
Lemma \ref{lem:faith:coef}, assume $a_1, \dots, a_{l'}$ are all
those maximal among $a_1,\dots,a_l$ and $b_{i,j}$ is maximal among
$\{ b_{i',j'} \mid 1 \leq i',j' \leq l' \}$. We claim $b_{i,i} =
b_{i,j} = b_{j,j}$.

The claim implies that $b_{i_1,\dots,i_m}$ is maximal among $\{
b_{i'_1,\dots,i'_m} \mid 1 \leq i'_1,\dots,i'_m \leq l' \}$ if and
only if $b_{i_j,i_k}$ is maximal among $\{ b_{i',j'} \mid 1 \leq
i',j' \leq l' \}$ for all $1 \leq j \leq k \leq l'$. Moreover, in
this case $\epsilon_{i_1} \cdots \epsilon_{i_m}
\rho_{n,m}(\alpha_{i_1,\dots,i_m})$ is independent of the choice
of $i_1,\dots,i_m$. Therefore, regarding $\pairrho{N}{\phi^{(m)}}$
as a polynomial of $q,t$, we find the coefficient of
$q^{a_{i_1,\dots,i_m}} t^{b_{i_1,\dots,i_m}}$ is nonvanishing thus
$\pairrho{N}{\phi^{(m)}} \neq 0$.

Now it remains to prove the claim. Let $\phi'$ denotes the other
component of $\phi^{(2)}$ and assume the tine edges of $\phi$ and
$\phi'$ intersect $N$ transversely at $\{ x_1,\dots,x_l \}$ and
$\{ x'_1,\dots,x'_l \}$, respectively. The rest part of the proof
is copied almost word by word from the proof of \cite[Claim
3.4]{Bigelow}.

Suppose, seeking a contradiction, that $b_{i,i}<b_{i,j}$. Let
$\alpha$ be an embedded arc from $z_i'$ to $z_j'$ along the tine
edge of $\phi'$. Let $\beta$ be an embedded arc from $z_j'$ to
$z_i'$ along $N$.

If $\beta$ does not pass through the point $z_i$, let
$\delta=\alpha\beta$ and let $w$ be the winding number of $\delta$
around $z_i$. Then $b_{i,j}-b_{i,i}=2w$. If $\beta$ does pass
through $z_i$, first modify $\beta$ in a small neighborhood of
$z_i$ so that $z_i$ lies to its left. Next let
$\delta=\alpha\beta$ and let $w$ be the winding number of $\delta$
around $z_i$. Then $1+b_{i,j}-b_{i,i}=2w$. In either case, our
assumption that $b_{i,i}<b_{i,j}$ implies that $w$ is greater than
zero.

Let $D_1=D\setminus\{z_i\}$. Let $\pi: \tilde D_1 \to D_1$ be the
universal (infinite cyclic) cover. Let $\tilde\alpha$ be a lift of
$\alpha$ to $\tilde D_1$. Let $\tilde\beta$ be the lift of $\beta$
to $\tilde D_1$ which starts at $\tilde\alpha(1)$. Let $\gamma$ be
a loop in $D_1$ based at $z_i'$ which winds $w$ times around $z_i$
in the clockwise (negative) direction such that $\gamma$ is
null-homotopic in $D \setminus P$. Let $\tilde\gamma$ be the lift
of $\gamma$ to an arc from $\tilde\beta(1)$ to $\tilde\alpha(0)$.
Choose $\gamma$ so that $\tilde\gamma$ is an embedded arc which
intersects $\tilde\alpha$ and $\tilde\beta$ only at its end
points.

Let $\tilde z_k'$ be the first point on $\tilde\alpha$ which
intersects $\tilde\beta$ (possibly $\tilde\alpha(1)$). Then
$\pi(\tilde z_k')=z_k'$ for some $k=1 \dots l$. Let
$\tilde\alpha'$ be the initial segment of $\tilde\alpha$ ending at
$\tilde z_k'$. Let $\tilde\beta'$ be the final segment of
$\tilde\beta$ starting at $\tilde z_k'$. Let
$\tilde\delta'=\tilde\alpha'\tilde\beta'\tilde\gamma$.

Now $\tilde\delta'$ is a simple closed curve in $\tilde D_1$, so
by the Jordan curve theorem it must bound a disk $\tilde B$. Since
$\gamma$ passes clockwise around $z_i$, there is a non-compact
region to the right of $\tilde\delta'$. Thus $\tilde\delta'$ must
pass counterclockwise around $\tilde B$.

Let $\alpha'$, $\beta'$ and $\delta'$ be the projections of
$\tilde\alpha'$, $\tilde\beta'$ and $\tilde\delta'$ to $D_1$. Then
$a_k-a_i$ is equal to the sum of the winding numbers of $\delta'$
around each of the points in $P$. This is equal to the cardinality
of $\tilde B \cap \pi^{-1}(P)$. Since $a_i$ is maximal among all
integers $a_{i'}$, we must have $a_k=a_i$. Thus $\tilde B \cap
\pi^{-1}(P)=\emptyset$. It follows that the arc
$\delta'=\alpha'\beta'\gamma$ is null-homotopic in $D \setminus
P$. But $\beta'$ is homotopic relative to end points to a subarc
of $N$, and $\gamma$ was chosen to be null-homotopic in $D
\setminus P$. Thus $\alpha'$ is homotopic relative to end points
to a subarc of $N$ in $D \setminus P$. So $\alpha$ and $N$ cobound
a digon in $D \setminus P$. But $\alpha'$ is a subarc of the tine
edge of $\phi$. This contradicts the fact that the tine edge of
$\phi$ intersects $N$ a minimal number of times. Therefore our
assumption that $b_{i,j}>b_{i,i}$ must have been false, so
$b_{i,j}=b_{i,i}$.

The proof that $b_{i,j}=b_{j,j}$ is similar. This completes the
proof of the claim, and hence of the lemma.
\end{proof}

Now we prove the main theorem.

\begin{figure}[h]
    \centering
    \includegraphics[scale=.6]{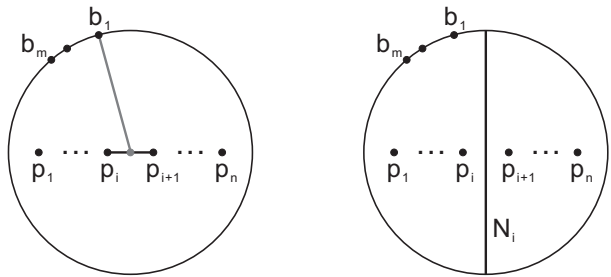}
    \caption{Fork $\phi_i$ and noodle $N_i$.}
    \label{fig:faith32}
\end{figure}

\begin{proof}[Proof of Theorem \ref{thm:faith}]
Suppose $\beta \in B_n$ belongs to the kernel of the Lawrence
representation, i.e. $L_{n,m}(\beta) = \id$. Then for any fork
$\phi$ and homeomorphism $f : D \to D$ representing $\beta$ , we
have $\pairrho{N}{\phi^{(m)}} = \pairrho{N}{(f\phi)^{(m)}}$.

Choose a set of disjoint noodles $N_1,\dots,N_{n-1}$ and a set of
forks $\phi_1,\dots,\phi_{n-1}$ with disjoint tine edges as shown
in Figure \ref{fig:faith32}. Note that
$\pairrho{N_j}{\phi_i^{(m)}} = 0$ if $i \neq j$. Choose a
homeomorphism $f$ representing $\beta$ such that $(f\phi_1)(e_t)
\cup \cdots \cup (f\phi_{n-1})(e_t)$ intersects $N_1 \cup \cdots
\cup N_{n-1}$ a minimal number of times in possible. Then,
whenever $i \neq j$, $\pairrho{N_j}{(f\phi_i)^{(m)}} = 0$ and by
Lemma \ref{lem:faith:key} $(f\phi_i)(e_t)$ is disjoint from $N_j$;
otherwise, $(f\phi_i)(e_t)$ and $N_j$ cobound a digon in $D
\setminus P$ which contradicts the minimality of the
intersections.

Therefore, we may further assume that $(f\phi_i)(e_t) =
\phi_i(e_t)$ thus $\beta$ must be a power of the full twist
$\Delta^2 = (\sigma_1 \cdots \sigma_{n-1})^n$. A straightforward
calculation shows that $L_{n,m}(\Delta^2) = q^{mn} t^{m(m-1)} \id$
hence we must have $\beta = 1$.
\end{proof}

\end{document}